\documentclass[12pt,reqno]{article}

\usepackage{amssymb,amsthm,amsmath,amsfonts,amstext}
\usepackage{latexsym,url}

\newtheorem{thm}{Theorem}
\newtheorem{cor}{Corollary}

\newtheorem{prop}{Proposition}

\newtheorem*{thmm}{Theorem}

\newcommand{\vv}{\vartheta}

\begin{document}

\begin{center}
\vskip 1cm{\LARGE\bf Congruence properties of the function which
counts compositions into powers of $2$ } \vskip 1cm \large Giedrius
Alkauskas\footnote{The author gratefully acknowledges support
from the Austrian Science Fund (FWF) under the project Nr. P20847-N18.}\\
Institute of Mathematics, Department of Integrative Biology\\
Universit\"{a}t f\"{u}r Bodenkultur\\
Gregor Mendel-Stra{\ss}e 33, A-1180 Wien, Austria, \&\\
Vilnius University, Department of Mathematics and
Informatics\\
Naugarduko 24, LT-03225 Vilnius, Lithuania.\\
{\tt giedrius.alkauskas@gmail.com} \\
\end{center}

\vskip .2in

\begin{abstract}
Let $\vv(n)$ denote the number of compositions (ordered partitions)
of a positive integer $n$ into powers of $2$. It appears that the
function $\vv(n)$ satisfies many congruences modulo $2^{N}$. For
example, for every integer $a$ there exists (as $k$ tends to
infinity) the limit of $\vv(2^k+a)$ in the $2-$adic topology. The
parity of $\vv(n)$ obeys a simple rule. In this paper we extend this
result to higher powers of $2$. In particular, we prove that for
each positive integer $N$ there exists a finite table which lists
all the possible cases of this sequence modulo $2^{N}$. One of our
main results claims that $\vv(n)$ is divisible by $2^{N}$ for almost
all $n$, however large the value of $N$ is.
\end{abstract}

\section{Introduction}

Let $\vv(n)$ denote the number of \emph{compositions} of a positive integer
$n$ into powers of $2$ (compositions are sometimes called \emph{ordered partitions}):
 this is the amount of finite sequences $\{q_{1},q_{2},\ldots,q_{\ell}\}$ of non-negative
 integers such that $n=2^{q_{1}}+2^{q_{2}}+\cdots+2^{q_{\ell}}$. Thus, for example, $
 3=1+1+1=2+1=1+2$ give all possible compositions, hence $\vv(3)=3$. This sequence appears
 in Sloane's encyclopedia \cite{sloane} as A023359. Several properties of this sequence are also listed there. Let us call this function \emph{the binary composition function}. It is easy to see (with the help of calculus of residues) that
\begin{eqnarray*}
\vv(n)\sim \frac{c}{\rho^{n+1}},
\end{eqnarray*}
where $\rho$ is the unique zero of $f(x)=1-\sum_{k=0}^{\infty}x^{2^{k}}$ in the interval $(0,1)$, and $c=-\frac{1}{f'(\rho)}$. Nevertheless, in this note we are mainly concerned with $2-$adic rather than with real asymptotics.\\

\emph{The binary partition function} $b(n)$ (which counts the partitions of $n$ into non-negative powers of $2$ neglecting the order of the summands) was investigated by many authors, beginning with L. Euler (1750), and in the 20th
century A. Tanturi (1918), K. Mahler (1940) (who explored asymptotic behavior). See the sequences A018819 and A000123 in Sloane's Encyclopedia \cite{sloane}; one can find numerous references there. Congruence properties of $b(n)$ modulo powers
of $2$ were first observed by R. F. Churchhouse \cite{church} (the main congruence was given without a proof as a conjecture). This conjecture was later proved by H. Gupta \cite{gupta} and independently by {\O}. R{\o}dseth \cite{rodseth}. This result can also be found in Andrews' monograph \cite{andrews}. The paper by the author \cite{alkauskas1} gives another proof of this fact along with one possible generalization of this congruence. As an aside \cite{alkauskas1}, for every positive integer $s$ which is not divisible by $8$ there exists a finite algorithm to verify the fact that infinitely many terms of the sequence $b(n)$ are divisible by $s$. Calculations confirmed this for $2\leq s\leq 14$, $s\neq 8$ (as was noticed by Churchhouse himself, $b(n)$ is never divisible by $8$). Moreover, for every power of $2$ there exists a finite table which lists all the possible remainders of $b(n)$ modulo this power. For example, modulo $32$ one of the entries is \cite{alkauskas1}
\begin{eqnarray*}
b(4n+2)\equiv 2+4w(n)+8w\big{(}\lfloor n/2\rfloor\big{)}+16\tau(n)\text{ (mod }32).
\end{eqnarray*}
Here $\lfloor\star\rfloor$ stands for the ``floor" function, $w(n)$ represents the Thue-Morse sequence with initial conditions $w(0)=0$, $w(1)=1$, and $\tau(n)$ stands for the Rudin-Shapiro sequence with conditions $\tau(0)=0$ and $\tau(3)=1$.\\

We will now formulate the R{\o}dseth-Gupta theorem.

\begin{thmm} If $n$ is odd positive integer, then for any integer $s\geq 1$ we have
\begin{eqnarray*}
b(2^{s+2}n)\equiv b(2^{s}n)\text{ (mod }2^{\mu(s)}),\quad
\text{were }\mu(s)=\Big{\lfloor}\frac{3s+4}{2}\Big{\rfloor};
\end{eqnarray*}
moreover, this congruence is exact.
\end{thmm}

On the other hand, the binary composition function has not yet been arithmetically investigated. The only papers (apart from The On-Line Encyclopedia of Integer Sequences) where this sequence appears are the papers by the author \cite{alkauskas2} and by Chinn and Niederhausen \cite{chinn}. The authors of the latter are concerned with finding an exact formula for the number of binary compositions of $n$ into exactly $n-k$ parts for small $k$.\\
\indent If we consider compositions of a positive integer $n$ with no limitation on the non-negative summands,
then the amount of these is equal to $2^{n-1}$. On the other hand, compositions with summands coming from a certain set reveal new congruence phenomena. For example, one of our main results is the following surprising fact. Let us denote by $s_{2}(n)$ the amount of $1$'s in the binary expansion of $n$. This is the sequence A000120.
\begin{thm}Suppose that $n\geq 1$, $N\geq 1$, and $s_{2}(n+2^{N-1})\geq 2^{N}$. Then
\begin{eqnarray*}
\vv(n)\equiv0\text{ (mod }2^{N}).
\end{eqnarray*}
\label{teom}
\end{thm}
\indent Let us say that a property $\mathcal{A}$ is satisfied for \emph{almost all natural numbers}, if
\begin{eqnarray*}
\lim\limits_{M\rightarrow\infty}\frac{\#\{n\leq M: n\text{ satisfies property }\mathcal{A}\}}{M}=1.
\end{eqnarray*}
\begin{cor}
Let $N\in\mathbb{N}$. Then for almost all natural numbers the congruence $\vv(n)\equiv0\text{ (mod }2^{N})$ is satisfied.
\end{cor}
\begin{proof} This is clear: for a fixed $M\in\mathbb{N}$, almost all natural numbers will have more than $M$ $1$'s in the binary expansion. \end{proof}
\section{Congruence properties}
Now we will derive some basic facts about $\vv(n)$. Les us make a convention $\vv(0)=1$ and $\vv(-n)=0$ for $n\in\mathbb{N}$.
Binary compositions of $n$ can be
divided into disjoined subsets, each of which consists of compositions with the first
summand equal to $2^k$, $1\leq 2^{k}\leq n$. Thus, this gives the recurrence relation, which also appears in \cite{alkauskas2,chinn,sloane}:
\begin{eqnarray}
\vv(n)=\sum_{k\geq0}\vv(n-2^k).\label{rec}
\end{eqnarray}
Hence, the generating function is give by
\begin{eqnarray*}
\sum_{n=0}^{\infty}\vv(n)x^{n}=(1-\sum_{k=0}^{\infty}x^{2^{k}})^{-1}.
\end{eqnarray*}
From the recurrence relation we can already determine the parity
of $\vv(n)$. This is the only property which admits easy proof directly from (\ref{rec}).
\begin{prop}
The number $\vv(n)$ is odd if and only if $n=2^{u}-1$,
$u\geq0$.\label{elem}
\end{prop}
\begin{proof} Suppose we have already proved this statement for all positive integers $\leq n-1$.
From the recurrence relation and inductive hypothesis we inherit that the parity of $\vv(n)$
equals the parity of the amount of odd terms among $\vv(n-2^{u})$, $u\geq 0$. This term
is odd iff (according to the induction hypothesis) $n-2^{u}=2^{v}-1$; that is, it happens only iff $n+1=2^{u}+2^{v}$.
Hence, if $s_{2}(n+1)>2$, this
cannot occur. If $s_{2}(n+1)=2$, so $n+1=2^{u}+2^{v}$, $u\neq v$, we have
exactly two odd summands: $\vv(n-2^{u})$ and $\vv(n-2^{v})$, and therefore the sum is even.
Finally, we have only one odd summand iff $n=2^{u}-1$, and this summand is
$\vv(n-2^{u-1})=\vv(2^{u-1}-1)$. We finish by induction. \end{proof}

\indent The congruence properties of the binary composition function modulo higher powers of $2$ were observed  by the author \cite{alkauskas2}. One of these congruences claim that
\begin{eqnarray}
\vv(2^{k})\equiv 8\text{ (mod }16)\text{ for }k\geq 3.\label{ast}
\end{eqnarray}
Unfortunately, despite many efforts to manipulate with (\ref{rec}), this and similar claims were not proved in \cite{alkauskas2} but rather extrapolated from numerical data. This failure suggests the fact that the recurrence (\ref{rec}) alone is insufficient in proving these congruences. Luckily, one can derive some other recurrence relations, much more convenient and powerful.
\begin{prop} For $n\geq1$, we have
\begin{eqnarray}
\vv(2n)&=&\vv^{2}(n)+\sum\limits_{{a+b=2n-2^{s}}\atop{a,b<n,s\geq 1}}\vv(a)\cdot\vv(b).\label{dvig}
\end{eqnarray}
In general, for $m,n\geq 1$, the following equality holds
\begin{eqnarray}
\vv(m+n)&=&\vv(m)\cdot\vv(n)+\sum\limits_{{a+b=m+n-2^{s}}\atop{a<m,b<n,s\geq 1}}\vv(a)\cdot\vv(b).\label{sum}
\end{eqnarray}
Thus, in case $m=n$ this gives (\ref{dvig}), and in case $m=1$ this reduces exactly to (\ref{rec}).
\end{prop}
\begin{proof}As a matter of fact, this identity is valid for any function which counts compositions of $n$ into positive integers $1=a_{1}<a_{2}<a_{3}<\cdots$, only $2^{s}$ must be replaced with $a_{s}$ in the formula. To prove the identity, consider any composition of $m+n$:
\begin{eqnarray*}
m+n=2^{q_{1}}+2^{q_{2}}+\cdots+2^{q_{\ell}}.
\end{eqnarray*}
Let $s$ be the largest non-negative integer such that $\sum_{i=1}^{s}2^{q_{i}}\leq m$ (if $s=0$, the empty sum is $0$ by convention). The number of compositions of $m+n$ where $\sum_{i=1}^{s}2^{q_{i}}=m$ for some $s$ is obviously equal to $\vv(m)\cdot\vv(n)$. If $a=\sum_{i=1}^{s}2^{q_{i}}<m$, then $b=\sum_{i=s+2}^{\ell}2^{q_{i}}<n$. Thus, $m+n-a-b=2^{q_{s+1}}$. Thus, if $a<m$ and $b<n$ are fixed, the amount of such compositions is equal to $\vv(a)\cdot\vv(b)$. This proves the formula (\ref{sum}). \end{proof}

Now we are able to derive the following
\begin{prop} The sequence $\vv(n)$ can be completely described modulo $4$.\\
(i) For $n\geq 3$, let $\tau_{3}(n)$ denote the number of solutions of $n+1=2^{s}+2^{v}+2^{u}$ with $s\geq v>u\geq 0$, let $\tau_{2}(n)$ denote the number of solutions of $n+1=2^{s}+2^{u}$, and let $\tau_{1}(n)$ denote the number of solutions of $n+1=2^{s}$. Then
\begin{eqnarray*}
\vv(2n)\equiv2\tau_{3}(n)+\tau_{2}(n)+\tau_{1}(n)\text{ (mod }4).
\end{eqnarray*}
(ii) In a similar way, $\vv(2^{k}-1)\equiv3\text{ (mod }4)$ for $k\geq 2$, and  $\vv(2^{k}+2^{l}-1)\equiv2\text{ (mod }4)$ for $k>l\geq1$. In all other cases, $\vv(2n-1)\equiv0\text{ (mod }4)$ (the first occurrence is $n=7$).
\label{proposition3}
\end{prop}
\begin{proof}{\it (i) }Note that (\ref{dvig}) can be rewritten as
\begin{eqnarray}
\vv(2n)&=&\vv^{2}(n)+2\sum\limits_{{a+b=2n-2^{s}}\atop{0\leq a<b<n,s\geq 1}}\vv(a)\cdot\vv(b)+\sum\limits_{s\geq 1}\vv^{2}(n-2^{s-1}).\label{even}
\end{eqnarray}
Consider this equality modulo $4$. Obviously, $\vv^{2}(n)\equiv \tau_{1}(n)\text{ (mod }4)$, since this is implied by Proposition \ref{elem}. In the second sum, only the terms with $a=2^{u}-1$ and $b=2^{v}-1$, $u,v\geq 0$, do contribute to the final result.
In this case $2n=2^{s}+2^{u}-1+2^{v}-1$. Thus, suppose $n$ is of this form. Since $a<b$, we have $u<v$, $u\geq 1$, and also, since $2^{v}-1<n$, it is easy to see that $s\geq v$. The number of solutions is thus $\tau_{2}(n)$. Finally, the second sum contributes exactly $\tau_{2}(n)$.\\

\noindent{\it (ii) }Equally, (\ref{sum}) for $m=n-1$ reads as
\begin{eqnarray}
\vv(2n-1)&=&\vv(n-1)\cdot\vv(n)+\sum\limits_{{a+b=2n-1-2^{s}}\atop{a<n-1,b<n,s\geq 1}}\vv(a)\cdot\vv(b)\nonumber\\
&=&\vv(n-1)\cdot\vv(n)+\sum\limits_{s\geq 1}\vv(n-2^{s})\cdot\vv(n-1)+2\sum\limits_{{a+b=2n-1-2^{s}}\atop{a<b<n-1,s\geq 1}}\vv(a)\cdot\vv(b)\nonumber\\
&=&2\vv(n-1)\cdot\vv(n)-\vv^2(n-1)+2\sum\limits_{{a+b=2n-1-2^{s}}\atop{a<b<n-1,s\geq 1}}\vv(a)\cdot\vv(b).\label{odd}
\end{eqnarray}
Here we used (\ref{rec}). In the sum, we have zero contribution to $\vv(2n-1)$ modulo $4$ unless $a=2^{u}-1$, $b=2^{v}-1$, and $2^{u}-1+2^{v}-1=2n-1-2^{s}$. Thus, $u=0$ and $n=2^{s-1}+2^{v-1}$. Since $b=2^{v}-1<n-1$, this implies $s>v$. Thus, there exists at most one such solution, and this happens exactly when $s_{2}(n)=2$. If $n=2^{k}$, $k\geq 1$, we get that $\vv(2^{k+1}-1)\equiv-\vv^{2}(2^{k}-1)\equiv 3\text{ (mod }4)$. \end{proof}

The following table summarizes the results of Proposition \ref{proposition3}.\\
\begin{center}
\begin{tabular}{|r | r|r|}
\hline
\multicolumn{3}{|c|}{\textbf{1. Sequence $\vv(n)$ modulo $4$, $n\geq 2$}}\\
\hline
$n$ & $\vv(n)\text{ (mod }4)$ & $\text{Condition}$\\
\hline
$2^{k}+2^{l}+2^{m}-2$ & $2$& $k>l>m\geq 1$ \\
$2^{k}-2$ & $2$& $k\geq 2$\\
$3\cdot2^{k}-2$ & $2$& $k\geq 1$\\
$\text{Other even numbers}$ & $0$& $ $\\
$2^{k}-1$ & $3$ & $k\geq 2$\\
$2^{k}+2^{l}-1$ & $2$ & $k>l\geq 1$\\
$\text{Other odd numbers}$ & $0$ & $ $\\
\hline
\end{tabular}\\
\end{center}

This table should list all even numbers $n$ such that $s_{2}(n+2)\leq 3$. However, two types of these numbers, namely $\{2^{k},k\geq 3\}$ and $\{2^{k}+2^{l}-2,k>l+1\geq 3\}$, fall under the qualification ``other even numbers", while the type $\{n=2^{k}+2^{l},\l\geq 2\}$ is a special case of the first type listed in the table.\\

Let us now inspect the recurrences (\ref{even}) and (\ref{odd}) more carefully. We will use the following well-known implication which, as a matter of fact, makes the investigations of quadratic forms over $2-$adic number field rather exceptional in $p-$adic analysis. Let $U,V\in\mathbb{Z}$, $N\geq 1$. Then
\begin{eqnarray*}
\text{if }U\equiv V\text{ (mod }2^{N}) \Rightarrow U^{2}\equiv V^{2}\text{ (mod }2^{N+1}).
\end{eqnarray*}
Suppose we know the sequence $\vv(n)$ modulo $2^{N}$. In this case the recurrences (\ref{even}), (\ref{odd}) and the above fact show us that the sequence $\vv(n)$ is completely describable modulo $2^{N+1}$ as well. Further, note that Table 1 lists only those even and odd numbers $n$ such that $s_{2}(n+2)\leq3$. The recurrence (\ref{even}) shows then that the corresponding table for $\vv(n)$ modulo $8$ will list only those even numbers $n$ such that $s_{2}(n+4)\leq 7$. Exactly the same conclusion follows for odd $n$. Here is one tricky point. In fact, consider odd number $2n-1$ and the multiplier $\vv(n)$ of the term $2\vv(n-1)\vv(n)$ in (\ref{odd}). This multiplier matters if $s_{2}(n+2)\leq 3$. This shows that odd numbers $2n-1$ such that $s_{2}((2n-1)+5)\leq 3$ should also be considered as candidates to be listed in the table for $\vv(n)$ modulo $8$. Indeed, it can happen that $s_{2}((2n-1)+5)\leq 3$ and $s_{2}((2n-1)+4)\geq 8$ are simultaneously satisfied. But then elementary considerations show that $s_{2}(n+1)\geq 7$. Thus, $\vv(n-1)\equiv 0\text{ (mod }4)$ and, due to this multiplier, the term $2\vv(n-1)\vv(n)$ does not contribute to $\vv(2n-1)\text{ (mod }8)$. We can proceed by induction on $N$. Therefore, a careful analysis of (\ref{even}) and (\ref{odd}) implies the following
\begin{thm}
For each positive integer $N$ there exists a finite table (analogous to the Table 1) which lists a finite number of possibilities for $\vv(n)$ modulo $2^{N}$. The table encompasses only finite number of classes of those positive integers $n$ such that
\begin{eqnarray}
s_{2}(n+2^{N-1})<2^{N}.
\label{bound}
\end{eqnarray}
 If the entry $n=2^{k_{1}}+2^{k_{2}}+\cdots+2^{k_{\ell}}-2^{N-1}$ is in this table, the corresponding residue depends solely  on $\ell$ and the exact shape of the collection of inequalities or equalities (the amount of these collections is also finite) satisfied by $k_{1},k_{2},\ldots,k_{\ell}$. These inequalities or equalities are of the form $k_{i}=k_{i+1}+d_{i}$, or  $k_{i}\geq k_{i+1}+d_{i}$, for a fixed collection of $d_{i}\in\mathbb{N}$. For those positive integers $n$ which are not in this table, $\vv(n)\equiv0\text{ (mod }2^{N})$.
\label{thm3}
\end{thm}
This result has numerous corollaries. One immediate corollary is Theorem \ref{teom}. Also, Theorem \ref{thm3} says that to prove the  congruence (\ref{ast}) or even to improve it to
\begin{eqnarray*}
\vv(2^{k})\equiv 8\text{ (mod }32)\text{ for }k\geq 8
\end{eqnarray*}
(which does hold), one needs to perform only a finite number of calculations: all what is demanded is to check that this congruence holds for $k$ up to a given bound, to be rest assured that this holds throughout. Indeed, according to the Theorem \ref{thm3}, two numbers $2^{k+1}$ and $2^{k}$ will eventually qualify for the same entry in the table for $k$ large enough (or be both left out of the table). For the very same reason, this allows to make the following claim. Let $a\in\mathbb{Z}$. Then there exists
\begin{eqnarray*}
\lim\limits_{k\rightarrow\infty}\vv(2^{k}+a)=\Theta(a)\in\mathbb{Z}_{2};
\end{eqnarray*}
here $\mathbb{Z}_{2}$ stands for the ring of $2-$adic integers, and the limit is taken in $2-$adic topology. A generalization of this is the following
\begin{cor}
Let $P(x)=\sum_{i=0}^{d}a_{i}x^{i}$ be a polynomial with non-negative integral coefficients. Then for every integer $N\in\mathbb{N}$ there exists $k_{0}\in\mathbb{N}$ such that
\begin{eqnarray*}
\vv\Big{(}P(2^{k+1})\Big{)}\equiv \vv\Big{(}P(2^{k})\Big{)}\text{ (mod }2^{N})\text{ for }k\geq k_{0}.
\end{eqnarray*}
\end{cor}
According to Theorem \ref{teom}, if a polynomial $P(x)=\sum_{i=0}^{d}a_{i}x^{i}$, $a_{d}\geq 1$, has at least one negative coefficient $a_{i}<0$ with $i\geq 1$, then $\vv(P(2^{k}))\rightarrow0$ in $2-$adic topology.\\

We finish with providing the table for $\vv(n)$ modulo $8$.

\begin{center}
\begin{tabular}{|r | r|r|}
\hline
\multicolumn{3}{|c|}{\textbf{2. Sequence $\vv(n)$ modulo $8$, $n\geq 7$}}\\
\hline
$n$ & $\vv(n)\text{ (mod }8)$ & $\text{Condition}$\\
\hline
$2^{k}+2^{l}+2^{m}-2$ & $6$& $k>l>m\geq 1$ \\
$2^{k}-2$ & $2$& $k\geq 2$\\
$3\cdot2^{k}-2$ & $6$& $k\geq 1$\\
$\text{Other even numbers}$ & $0$& $ $\\
$2^{k}-1$ & $7$ & $k\geq 3$\\
$2^{k}+1$ & $6$ & $k\geq 4$\\
$2^{k}+2^{l}-1$ & $2$ & $k>l\geq 2$\\
$2^{k}-3$ & $4$ & $k\geq 4$\\
$3\cdot2^{k}-3$ & $4$ & $k\geq 3$\\
$2^{k}+2^{l}+2^{m}-3$ & $4$ & $k>l>m\geq 2$\\
$\text{Other odd numbers}$ & $0$ & $ $\\
\hline
\end{tabular}\\
\end{center}
Thus, for example, $4|\vv(2n)\Rightarrow8|\vv(2n)$. The complete table for $2^{N}=16$ has the following entries, which are not members of larger classes:
\begin{eqnarray*}
\vv(7\cdot2^{k}-2)\equiv 14 \text{ (mod }16)\text{ for }k\geq 1, \quad\vv(5\cdot2^{k}-2)\equiv 8 \text{ (mod }16)\text{ for }k\geq 3.
\end{eqnarray*}
One could wonder, for example, whether these have combinatoric proofs.\\

\noindent Concerning the condition (\ref{bound}), it is certainly not sharp for $N>2$, as the above table for $2^{N}=8$ suggests. Employing this table, we can show that in fact $2^{N}$ on the right hand side can be replaced by $2^{N-1}+2^{N-3}-1$ for $N\geq 4$. This might be also far from the optimal for larger $N$.

\noindent\begin{center}
\begin{tabular}{|r| r r||r| r r||r| r r|}
\hline
\multicolumn{9}{|c|}{\textbf{4. Sequence $\vv(n)$}}\\
\hline
$n$ & $\vv(n)$& mod $64$& $n$ & $\vv(n)$& mod $64$ & $n$ & $\vv(n)$& mod $64$\\
\hline
$1$ & $1$&$000001$        &$25$&$ 882468$&$100100$ &$49$&$751322695068$&$011100$\\
$2$ & $2$&$000010$        &$26$&$ 1558798$&$001110$ &$50$&$1327134992166$&$100110$\\
$3$ & $3$&$000011$        &$27$&$ 2753448$&$101000$ &$51$&$2344248747712$&$000000$\\
$4$ & $6$&$000110$        &$28$&$ 4863696$&$010000$ &$52$&$4140876568224$&$100000$\\
$5$ & $10$&$001010$       &$29$&$ 8591212$&$101100$ &$53$&$7314436562436$&$000100$\\
$6$ & $18$&$010010$       &$30$&$ 15175514$&$011010$ &$54$&$12920206953182$&$011110$\\
$7$ & $31$&$011111$       &$31$&$ 26805983$&$011111$ &$55$&$22822229201360$&$010000$\\
$8$ & $56$&$111000$       &$32$&$ 47350056$&$101000$ &$56$&$40313142631496$&$001000$\\
$9$ & $98$&$100010$       &$33$&$83639030$&$110110$ &$57$&$71209059135432$&$001000$\\
$10$ & $174$&$101110$     &$34$&$147739848$&$001000$ &$58$&$125783547796216$&$111000$\\
$11$ & $306$&$110010$     &$35$&$260967362$&$000010$ &$59$&$222183821668104$&$001000$\\
$12$ & $542$&$011110$     &$36$&$460972286$&$111110$ &$60$&$392465083678728$&$001000$\\
$13$ & $956$&$111100$     &$37$&$814260544$&$000000$ &$61$&$693249583836156$&$111100$\\
$14$ & $1690$&$011010$    &$38$&$1438308328$&$101000$ &$62$&$1224554757801706$&$101010$\\
$15$ & $2983$&$100111$    &$39$&$2540625074$&$110010$ &$63$&$2163051215343439$&$001111$\\
$16$ & $5272$&$011000$    &$40$&$4487755390$&$111110$ &$64$&$3820809588459176$&$101000$\\
$17$ & $9310$&$011110$    &$41$&$7927162604$&$101100$ &$65$&$6749070853108302$&$001110$\\
$18$ & $16448$&$000000$   &$42$&$14002525142$&$010110$ &$66$&$11921546029897416$&$001000$\\
$19$ & $29050$&$111010$   &$43$&$24734033936$&$010000$ &$67$&$21058196429732338$&$110010$\\
$20$ & $51318$&$110110$   &$44$&$43690150992$&$010000$ &$68$&$37197158469308174$&$001110$\\
$21$ & $90644$&$010100$   &$45$&$77174200244$&$110100$ &$69$&$65704990586807960$&$011000$\\
$22$ & $160118$&$110110$  &$46$&$136320361910$&$110110$ &$70$&$116061171489076784$&$110000$\\
$23$ & $282826$&$001010$  &$47$&$240796030130$&$110010$ &$71$&$205010234490786986$&$101010$\\
$24$ & $499590$&$000110$  &$48$&$425341653750$&$110110$ &$72$&$362129691668018062$&$001110$\\
\hline
\end{tabular}
\end{center}
\bigskip
\hrule
\bigskip

\noindent 2010 {\it Mathematics Subject Classification}: Primary
11P83; Secondary 11P81, 05A17

\noindent \emph{Keywords: } Binary compositions, ordered partitions,
congruence properties, $2$-adic analysis

\bigskip
\hrule
\bigskip

\noindent (Concerned with sequence  A023359)

\bigskip
\hrule
\bigskip

\end{document}